\documentclass{elsart}

\usepackage{amssymb,cite,graphicx,here}

\newtheorem{theo}{Theorem} 

\begin{document}

\begin{frontmatter}



\title{On numerical solutions to stochastic Volterra equations\thanksref{iccam}}

 \author[AK]{Anna Karczewska\corauthref{cor1}}
 \address[AK]{Department of Mathematics, University of Zielona G\'ora, 
 65-246 Zielona G\'ora, Poland}
\ead{~A.Karczewska@im.uz.zgora.pl}
\corauth[cor1]{Corresponding author}
\thanks[iccam]{Extended version of the talk given at {\sc International Congress on 
Computational and Applied Mathematics}, Leuven, July 26-30, 2004}

\author[PR]{Piotr Rozmej}
\ead{\\~~P.Rozmej@if.uz.zgora.pl}
\ead[url]{www.uz.zgora.pl/$\sim$prozmej}

\address[PR]{Institute of Physics, University of Zielona G\'ora, 65-246
 Zielona G\'ora, Poland}

\begin{abstract}
 The aim of the paper is to demonstrate the use of the Galerkin method for some
 kind of Volterra equations, determininistic and stochastic as well.
 The paper consists of two parts: the theoretical and numerical one.
 In the first part we recall some apparently well-known 
 results concerning the Volterra
 equations under consideration. In the second one we describe a numerical algorithm
 used and next present some examples of numerical solutions in order to illustrate 
 the pertinent features  of the technique used in the paper. 
\end{abstract}

\end{frontmatter}

\section{Introductiom}\label{s1}
In the paper we investigate a stochastic version of a linear Volterra equation 
of the general form
\begin{equation}\label{eq1}
X(t,x)=\int\limits_0^t a(t-\tau)AX(\tau,x)d\tau+ X_0(x)+ f(t,x),
\end{equation}
where $t\in\mathbb{R}_+$, $x\in \mathbb{R}^d$, $a\in L^1_{loc}(
\mathbb{R}_+)$, $A$ is a linear operator and $f$ some mapping.
To fix our attention we shall consider the equation (\ref{eq1}) in a separable Hilbert
space $H$ with a scalar product $(\cdot ,\cdot)$, a norm $|\cdot|$ and a complete
orthonormal system $\{e_n\}$. The equation (\ref{eq1}) creates a big class of equations
and generalizes heat and wave equations and even linear Navier-Stokes equation.
We refer to the excellent monograph \cite{Pr} for a rich survey.
That kind of Volterra equation has been studied by many authors in connection with
problems arising in mathematical physics, particularly in viscoelasticity,
heat conduction in materials with memory, energy balance and termoviscoelasticity.
In order to take into account random fluctuations, we have to consider the 
equation (\ref{eq1}) with random external force.

There are our first considerations concerning numerical treatment 
of stochastic Volterra equations, so we will be grateful for readers'
remarks and advices.

Next, we plan to study the probabilistic features of family of 
trajectories, take into account different noises and develope numerical schemes
for cases where the analytic form of reselvent is not known.

\section{Resolvent approach}

Assume that $(\Omega,\mathcal{F},(\mathcal{F}_t)_{t\geq 0},P)$ is a probability
space with a complete right-continuous filtration and $W(t),~t\geq 0$, is a
cylindrical Wiener process with values in the space $H$. Let us omit, for convenience,
the space variable $x$ in the equation (\ref{eq1}) and introduce the process 
$W(t),~t\geq 0$, instead of function~$f$. Hence, we arrive at the following 
stochastic Volterra equation
\begin{equation}\label{eq2}
X(t)=\int\limits_0^t a(t-\tau)AX(\tau)d\tau+ X_0+ W(t).
\end{equation}

In this part of the paper we recall some results concerning solutions to 
(\ref{eq2}). We restrict our considerations to paper containing so-called resolvent
approach to Volterra equation. The notion {\em resolvent} or {\em fundamental solution}
for Volterra equation (\ref{eq1}) probably comes from Friedman and Shinbrot \cite{FrSh}
who studied deterministic Volterra integral equations in Banach space. 
For recent survey we refer again to \cite{Pr}.
In the sequel we shall assume that the equation (\ref{eq1}) is 
\textit{well-posed}, that is, 
that (\ref{eq1}) admits resolvent $S(t),~t\geq 0$.

As in deterministic case, the {\it mild solution} 
to the stochastic Volterra equation (%
\ref{eq2}) is of the form
\begin{equation}  \label{eq3}
X(t) = S(t)X_0 + \int_0^t S(t-\tau)dW(\tau)\;, \quad
t\geq 0\;,
\end{equation}
where $S(t),~t\geq 0$, is the resolvent family for the equation 
(\ref{eq1}) determined by the operator $A$ and the function $a$.
In order to study solution (\ref{eq3}) it is enough to consider the
stochastic convolution 
\begin{equation}  \label{eq4}
W_{S}(t) := \int_0^t S(t-\tau)dW(\tau)\;,\quad t\geq 0\;,
\end{equation}
where the stochastic integral is defined according to particular case under
consideration.

Stochastic Volterra equations with resolvent approach have been treated by 
several authors, see e.g. 
\cite{ClDaPr1},\cite{ClDaPr2},\cite{ClDaPr3},\cite{ClDaPrPr},\cite{KaZa} 
and recently \cite{Ka1} and \cite{Ka2}. In the first
three papers stochastic Volterra equations are studied in connection with
viscoelasticity and heat conduction in materials with memory.
The paper  \cite{ClDaPr1} is particularly significant
because the authors were the first who 
have extended the well-known semigroup approach, applied
to stochastic differential equations, to the
equation (\ref{eq2}). The resolvent approach is a natural way of extension the 
semigroup approach which is well-known from the theory of evolution equations.
That approach enables to follow some results and schemes obtained for semigroups.
Unfortunately, some results are not valid in our case because resolvent family,
$S(t),~t\geq 0$, does not satisfy semigroup property.

Cl\'ement and DaPrato studied stochastic Volterra equation (\ref{eq2}) 
where $A$ was self-adjoint, 
negative operator in the space $H$, such that
\[
A e_k = - \mu_k e_k, \quad \mu_k>0, \quad k\in \mathbb{N}\;. 
\]
They considered stochastic Volterra
equation (\ref{eq2}) driven by the noise term $W$ of the form 
\begin{equation}  \label{eq5}
(W(t),h)_H = \sum_{k=1}^{+\infty} (h,e_k)_H\,\beta_k(t),\quad h\in H,
\end{equation}
where $\{\beta_k\}$ was a sequence of real-valued, independent Wiener
processes. They assumed that the kernel function $a$ is completely
positive. The consequence of completely positiveness of the function $a$ is
that the solution $s(\cdot,\gamma),\hspace{1ex}\gamma>0$, to the following
equation 
\begin{equation}  \label{eq6}
s(t) + \gamma\int_0^t v(t-\tau) s(\tau) d\tau = 1, \quad t\geq 0,
\end{equation}
is nonnegative and nonincreasing for any $\gamma>0$. In fact, 
$s(t)\in [0,1]$. 
In \cite{ClDaPr1} regularity of stochastic convolution (\ref{eq4}) 
is studied and h\"olderianity of the
corresponding trajectories is proved.

\begin{hypo} \label{hyp1}  
~\\ \vspace{-5mm}
 \begin{itemize}
  \item[(i)] A is a self-adjoint negative operator and $Ae_k=-\mu_k e_k$.
  \item[(ii)] $a$ is completely positive.
  \item[(iii)]$\quad -\mathrm{Tr}(A^{-1}) = \sum_{k=1}^{\infty}
       (1/\mu_k) < +\infty $.
 \end{itemize}
\end{hypo}

\begin{hypo} \label{hyp2} 
 There exists $\theta\in (0,1)$ and $C_\theta>0$ such that, \\
 for all $0<\tau<t$ we have
 $$ \int_{\tau}^t s^2(\mu ,\sigma)d\sigma \leq C_\theta 
    \mu^{\theta-1} |t-\tau|\theta\;,$$
 $$ \int^{\tau}_0 [s(\mu ,\tau-\sigma)-s(\mu ,t-\sigma)]^2 d\sigma
    \leq C_\theta\mu^{\theta} |t-\tau|\theta$$
and
 $$ \sum_{k=1}^{\infty} \mu_k^{\theta-1} < +\infty\;.$$
\end{hypo}
\begin{hypo} \label{hyp3} 
 There exists $M>0$ such that
 $$\left\{ \begin{array}{l}
  |e_k(\theta)| \leq M, \quad k\in \mathbb{N}, \quad \theta\in\mathcal{O}, \\
  |\nabla e_k(\theta)| \leq M\mu_k^{1/2}, \quad k\in \mathbb{N},
  \quad \theta\in\mathcal{O},
 \end{array} \right.$$
where $\mathcal{O}$ is a bounded open subset of $\mathbb{R}^d$.
\end{hypo}
Cl\'ement and DaPrato proved the following  results.
\begin{theo} \label{3.1} 
(\cite{ClDaPr1}, Theorem 2.2)
 Assume that Hypothesis \ref{hyp1} holds. Then for any $t\geq 0$
 the series
 $$ \sum_{k=1}^{\infty} \int_0^t s(\mu_k,t-\tau)e_k d\beta_k(\tau), $$
is convergent in $L^2(\Omega)$ to a Gaussian random variable $W_S(t)$ 
with mean 0 and covariance operator $Q_t$ determined by
$$ Q_t e_k = \int_0^t s^2(\mu_k,\tau)d\tau e_k, \quad k\in \mathbb{N}.$$
\end{theo}
\begin{theo} \label{3.2} 
(\cite{ClDaPr1}, Proposition 3.3)
 Under Hypotheses \ref{hyp1} and \ref{hyp2}, for every positive number
 $\alpha <\theta /2$, the trajectories of $W_S$ are almost surely
 $\alpha$-H\"older continuous. 
\end{theo}
\begin{theo} \label{3.3} 
(\cite{ClDaPr1}, Theorem 4.1)
 Under Hypotheses \ref{hyp1}, \ref{hyp2} and \ref{hyp3}, the trajectories 
 of $W_S$ are almost surely $\alpha$-H\"older continuous in
 $(t,x)$ for any $\alpha\in (0,1/4)$.
\end{theo}
In the paper \cite{ClDaPr2}, 
white noise perturbation of an integro-differential equation arising
in the study of evolution of material with memory is studied. 
In the next paper 
\cite{ClDaPrPr}, the authors considered evolutionary integral equations as
appearing in the theory of linear parabolic viscoelasticity forced by white
noise. As earlier, they studied the stochastic convolution that provides
regular solutions. Additionally, under suitable assumptions the authors
proved that the samples are H\"older-continuous. In the remaining part of
the paper \cite{ClDaPrPr}, the results obtained of that paper were put in a
wider perspective by consideration of equations with fractional derivatives.

In the paper \cite{ClDaPr3}, the authors
first proved that $W_S(t)$ is a Gaussian random variable for any $t\geq 0$.
Next, the transition function $P_t,\, t\geq 0$, associated with
$S(t)$ was considered. When $S(t)$ is the resolvent
operator of a stochastic Volterra equation, the convolution $W_S(t),\, t\geq
0$, is not a Markov process. This fact has the consequence that $P_t,\, t\geq
0$, is not a semigroup and then it is not possible to associate to  $P_t$ a
Kolmogorov equation. However, the authors characterized those
transition functions such that $P_t\,\varphi$ was differentiable for any
uniformly continuous and bounded function $\varphi$. 

There are some other regularity results concerning stochastic convolution (\ref{eq4}).
In \cite{KaZa} was studied the case when the equation (\ref{eq2}) was driven by a
correlated, spatially homogeneous Wiener process $W$ with values in the space of real,
tempered distributions $S'(\mathbb{R}^d)$. Let $\Gamma$ be the
covariance of $W(1)$ and the associated spectral measure be $\mu$.
We considered existence of the solutions to (\ref{eq2}) in 
$S^{\prime}(\mathbb{R}^d)$ and derived conditions under
which the solutions to (\ref{eq2}) were function-valued and continuous. In that
case,  the initial value $X_0 \in S'(\mathbb{R}^d)$, 
$a$ is a locally integrable function and $A$ is an operator
given in the Fourier transform form 
\begin{equation}\label{eq7}
\mathcal{F}(A\xi)(\lambda) = -\beta(\lambda)\mathcal{F}(\xi)(\lambda)\;.
\end{equation}
 We introduce the following hypothesis.
\begin{hypo} \label{H} 
{}
\begin{enumerate}
 \item For any $\gamma\geq 0$, the equation (\ref{eq6}) has exactly
  one solution $s(\cdot,\gamma)$ locally integrable and measurable
  with respect to both variables $\gamma\geq 0$ and $t\geq 0$.
 \item Moreover, for any $T\geq 0$, 
\quad
 $ \sup_{t\in [0,T]} \; \sup_{\gamma\geq 0} |s(t,\gamma)|<+\infty .$
\end{enumerate}  
\end{hypo}

\vspace{1mm} For some special cases the function
$s(t;\gamma)$ may be found explicitly. For instance 
\begin{equation}  \label{eq8}
\begin{array}{llll}
\mathrm{for} & a(t)\equiv 1, & s(t;\gamma)=e^{-\gamma t}, & ~t\geq 0,\hspace{%
1ex} \gamma\geq 0; \\ 
\mathrm{for} & a(t)=t, & s(t;\gamma)=\cos(\sqrt{\gamma}t), & ~t\geq 0,\hspace{%
1ex} \gamma\geq 0; \\ 
\mathrm{for} & a(t)=e^{-t}, & s(t;\gamma)=(1+\gamma)^{-1} [1+\gamma
e^{-(1+\gamma)t}], & ~t\geq 0,\hspace{1ex} \gamma\geq 0.
\end{array}
\end{equation}

 In that case, the resolvent family $S(t),~t\geq 0$,
determined by the operator $A$ and the function $a$ is given 
by the formula (\ref{eq7}) and has the form
\[ S(t)\xi = \mathbf{r}(t)\star \xi, \quad \xi\in S^{\prime}(%
\mathbb{R}^d), \]
where 
$\mathbf{r}(t) = \mathcal{F}^{-1} s(t,\beta(\cdot)),~t\geq 0. $
The following results for stochastic convolution are 
consequences of properties of stochastic integral.

\begin{theo} \label{3.4} 
(\cite{KaZa}, Theorem 2)
 Let $W$ be a spatially homogeneous Wiener process and 
 $S(t), t\geq 0$, the resolvent for the equation (\ref{eq2}).
 If Hypothesis \ref{H} holds then the stochastic equation 
 $$ S\star W(t) = \int_0^t S(t-\sigma)
 dW(\sigma), \quad t\geq 0$$
 is a well-defined $S'(\mathbb{R}^d)$-valued process. For each $t\geq 0$ the
 random variable $S\star W(t)$ is generalized, 
 stationary random field on $\mathbb{R}^d$ with the spectral measure $\mu_t$:
 $$\mu_t(d\lambda) = \left[ \int_0^t (s(\sigma,\beta(\lambda)))^2
 d\sigma \right] \mu(d\lambda)\;.$$
\end{theo}
\begin{theo} \label{3.5} 
(\cite{KaZa}, Theorem 3)
 Assume that Hypothesis \ref{H} holds. Then the process 
 $S\star W(t)$ is function-valued for all $t\geq 0$
 if and only if 
 $$ \int_{\mathbb{R}^d}\left( \int_0^t (s(\sigma,\beta(\lambda)))^2 d\sigma 
 \right) \mu(d\lambda)<+\infty , \quad  t\geq 0.$$
 If for some $\epsilon >0$ and all $t\geq 0$,
 $$ \int_0^t \int_{\mathbb{R}^d} (\ln (1+|\lambda|))^{1+\epsilon}
  (s(\sigma,\beta(\lambda)))^2 d\sigma \mu(d\lambda)<+\infty , $$
 then, for each $t\geq 0$, $S\star W(t)$ is a sample
 continuous random field.
\end{theo}

As we have already written, the Volterra equation (\ref{eq1}) creates
a big class of equations. In particular cases, when the
operator $A$ and the function $a$ in the stochastic Volterra equation 
(\ref{eq2}) are fixed, there is possible to obtain some additional
regularity results. This is obvious because in particular cases we may
use some extra features of solutions to (\ref{eq2}).
For instance, apparently well-known is integrodifferential equation
which interpolates heat and wave equations, that is the 
Volterra equation (\ref{eq1}), where $A=\Delta$, the Laplace operator
and $a(t)=t^{\alpha-1}/\Gamma(\alpha),~0\leq \alpha\leq 2$, where
$\Gamma$ is the gamma function. 
Recently, deterministic version was studied in detail by Fujita 
\cite{Fu} and, independently, by Schneider and Wyss \cite{SchWy}
and stochastic version of that integrodifferential equation 
was treated in \cite{Fu2} and \cite{Ka2}. In this paper 
we shall demonstrate numerical results obtained for that equation
in the deterministic version and the stochastic one, as well.

In the theory of stochastic Volterra equations we are interested not
only in the existance, uniqueness and regularity of solutions but 
in some asymptotics, too. The paper \cite{Ka1} is concerned with 
a limit measure of stochastic Volterra equation driven by very general
noise in a form of a spatially homogeneous Wiener process, with values
in the space of tempered distributions $S'(\mathbb{R}^d)$.
That paper provides necessary and sufficient conditions for the
existence of the limit measure and additionally, it gives a form of any
limit measure.

Let us summarize the results cited above. 
The resolvent operators $S(t),~t\geq 0$, corresponding to the
Volterra equation (\ref{eq1}) do not form any semigroup. Therefore it
is not possible to obtain such strong results as in the case of 
evolution equations with semigroup generators.
In the case of Volterra equations one can not use 
{\it the fractional method of infinite dimensional
stochastic calculus}. That method, used for demonstration of continuity 
with respect to $t$ for convolutions with semigroups, enables to obtain
only some estimates for the convolutions  (\ref{eq4}). Moreover, it is
possible only in some special cases, see e.g.\ \cite{Ka3}.
It is clear that the Volterra equation (\ref{eq1}), in particular its
stochastic version (\ref{eq2}), is difficult to study. It results from
the fact that equation (\ref{eq1}) contains a wide class of equations. 
An essential role is played by the kernel function 
which, in general, is assumed to be a locally
integrable function. As the function $a$ and operator $A$ determine the
resolvent $S(t),~t\geq 0$, the type of the function $a$ is very 
important. That is reason why it is so difficult to obtain in 
a general case the continuity of the convolution (\ref{eq4}) 
and some other theoretical results.
For more general convolution, significant
in many applications, like $W^\psi(t):=\int_0^t S(t-\tau)\psi(\tau)dW(\tau)$, 
where $\psi$ is an appropriate process, 
it becomes even far more difficult.

Therefore, in many cases a numerical support of 
theoretical (analytical) considerations is demanded. 
We need computations for obtaining estimates in regularity results, 
choosing some paramenters, choosing function $a$ with
required properties and for visualization of solutions obtained.
Numerical analysis is particularly important when we are not able to
obtain analytical results. Additionally, numerical schemes are 
especially useful for studying the asymptotics.

\section{Galerkin method for deterministic Volterra equation}
\label{s2}

In this section we construct a scheme for numerical solution of the 
Volterra  equation (\ref{eq2})
without random part, that is, for
 \begin{equation}\label{eq1a}
X(t,x)=\int\limits_0^t a(t-\tau)AX(\tau,x)d\tau+ X_0(x).
\end{equation}
We shall consider the case when $A$ is the Laplace operator.
Denoting by $K(x,t,s)= a(t-s)\frac{d^2}{dx^2}$ we can write (\ref{eq1a})
in the standard form
\begin{equation} \label{e2}
X(x,t) = X_0(x) + \int_0^t K(x,t,s)\,X(x,s) ds \;.
\end{equation}

In Galerkin method one introduces the complete set of orthonormal functions
$\{ \phi_j\},~~j =1,\ldots ,\infty$ on the interval  $[0,t]$, that is fulfilling
conditions
\begin{equation} \label{e4}
(\phi_i(t),\phi_j(t))= \int_0^t \phi_i(\tau)\,\phi_j(\tau)\,d\tau = \delta_{ij}  \;,
\end{equation}
where  $(\cdot,\cdot)$ is the scalar product. The set $\{ \phi_j\}$ spans 
a Hilbert space. The approximate solution is then postulated in the form
of an expansion of the unknown true solution in the subspace $H_n$ determined
by $n$ first basis functions
\begin{equation} \label{e3}
X_n(x,t) =\sum_{j=1}^n c_j(x)\,\phi_j(t) \;.
\end{equation}
Inserting (\ref{e3}) into (\ref{e2})  we obtain
\begin{equation} \label{e5}
X_n(x,t) = X_0(x) + \int_0^t  K(x,t,s)\,X_n(x,s) ds +\varepsilon_n(x,t) \;,
\end{equation}
where the function $\varepsilon_n(x,t)$ represents the approximation error. 
From (\ref{e5}) we have
\begin{equation} \label{e6}
\varepsilon_n(x,t)=f_n(x,t) - X_0(x) - \int_0^t  K(x,t,s)\,f_n(x,s) ds \;.
\end{equation}
From  (\ref{e6}) and (\ref{e3}) it can be written as
\begin{equation} \label{e7}
\varepsilon_n(x,\tau)=\sum_{k=1}^n c_k(x)\,\phi_k(\tau)  - X_0(x)
 - \int_0^\tau  K(x,\tau,s)\, \sum_{k=1}^n c_k(x)\,\phi_k(s)\,ds \;.
\end{equation}
The coefficient functions  $c_j(x)$ are determined by the requirement that
the error function $\varepsilon_n(x,t)$ has to be orthogonal to the subspace $H_n$
\begin{equation} \label{e8}
(\phi_j(t), \varepsilon_n(x,t)) =0\quad\quad j =1,2,\ldots,n \;.
\end{equation}
Then for $j =1,2,\ldots,n$ the following equations hold
\begin{eqnarray} \label{e9}
 \int_0^t X_0(x)\phi_j(\tau) d\tau  &=& 
 \int_0^t \left[\sum_{k=1}^n c_k(x)\,\phi_k(\tau)\right] \phi_j(\tau)d\tau
 \\ & -&  
 \int_0^t \left[ \int_0^\tau K(x,\tau,s) \sum_{k=1}^n c_k(x)\,\phi_k(s)ds
\right] \phi_j(\tau)d\tau\;. \nonumber
\end{eqnarray}
The first integral on the r.h.s, due to (\ref{e3}), is very simple
\begin{equation} \label{e10}
 \int_0^t \left[\sum_{k=1}^n c_k(x)\,\phi_k(\tau)\right]\phi_j(\tau)d\tau 
 = \sum_{k=1}^n c_k(x) \int_0^t \phi_k(\tau)\phi_j(\tau)d\tau 
  = c_j(x) \;.
\end{equation}
The second one we calculate in our particular case, 
$K(x,\tau,s)= a(\tau-s)\frac{d^2}{dx^2}$ as follows
\begin{eqnarray} \label{e11}
\int_0^t \left[ \int_0^\tau \sum_{k=1}^n \frac{d^2 c_k(x)}{dx^2}\,a(\tau-s)
 \phi_k(s)ds \right] \phi_j(\tau)d\tau \hspace{25ex} \\  \hspace{18ex}
 = \sum_{k=1}^n \frac{d^2 c_k(x)}{dx^2}
 \int_0^t \left[ \int_0^\tau a(\tau-s)\phi_k(s)ds \right] \phi_j(\tau)d\tau
 \nonumber \;.
\end{eqnarray}

Denoting by
\begin{eqnarray} \label{e11a}
 g_j(x) & = &  \int_0^t X_0(x)\phi_j(\tau) d\tau =
  X_0(x)\int_0^t \phi_j(\tau) d\tau\;, \quad \quad \mbox{and} \\
 a_{jk} & = &
\int_0^t \phi_j(\tau)\left[ \int_0^\tau a(\tau-s)\phi_k(s)ds \right] d\tau 
\end{eqnarray} 
we arrive at the set of coupled differential equations for the functions $c_j(x)$
\begin{equation} \label{e12}
 g_j(x) = c_j(x) - \sum_{k=1}^n  a_{jk}\frac{d^2 c_k(x)}{dx^2}  \;.
\end{equation}
This set can be solved numerically (aproximately) by discretization,
on a grid  $x_i=\{x_1,x_2,\ldots,x_m\}$. Applying the difference form for the
second derivative
$$ \frac{d^2 c_k(x_i)}{dx^2} \approx \frac{1}{h^2} 
  [c_k(x_{i-1})-2 c_k(x_i) + c_k(x_{i+1})], \quad\quad h= x_{i}-x_{i-1}
$$ 
one obtains from (\ref{e12}) the following set of linear equations
\begin{eqnarray} \label{e13}
 g_j(x_i) & =& c_j(x_i) + \frac{1}{h^2} \sum_{k=1}^n a_{jk}\, 
 [-c_k(x_{i-1})+2 c_k(x_i) - c_k(x_{i+1})] \;,
\end{eqnarray}
where $j=1,2,\ldots,n,~i=1,2,\ldots,m$. Those equations can be written in
the matrix form
\begin{equation} \label{e14}
 \underline{\mathcal{A}}\, \underline{c} = \underline{g}\;.
\end{equation}
Here,  $(N=n\times m)$-dimensional vectors 
$\underline{c}$ i $\underline{g}$ have the following structures
\begin{eqnarray} \label{e15}
 \underline{c} =\left( \begin{array}{c}
 c_1(x_1) \\ \vdots \\ c_1(x_m)\\  c_2(x_1) \\ \vdots \\ c_2(x_m)\\  \cdots \\
 c_n(x_1) \\ \vdots \\ c_n(x_m)\\ 
 \end{array}\right) = \left( \begin{array}{c} C_1 \\ C_2 \\ \vdots \\ C_n
  \end{array}\right)
 \;, \hspace{10ex}
  \underline{g} =\left( \begin{array}{c}
 g_1(x_1) \\ \vdots \\ g_1(x_m)\\  g_2(x_1) \\ \vdots \\ g_2(x_m)\\  \cdots \\
 g_n(x_1) \\ \vdots \\ g_n(x_m)\\ 
 \end{array}\right) = \left( \begin{array}{c} G_1 \\ G_2 \\ \vdots \\ G_n
  \end{array}\right) \; ,
\end{eqnarray}
where by $C_i$, $G_i$ we denoted the consecutive $m$-dimensional
blocks of vectors $\underline{c}$ and $\underline{g}$, respectively.
Then we can write the matrix $\underline{\mathcal{A}}$ in the block form
\begin{eqnarray} \label{e16}
\underline{\mathcal{A}} =  \left( \begin{array}{ccc}
 \left[A_{11}\right] & \ldots & \left[A_{1n}\right] \\ 
 \vdots & \cdots & \vdots \\
 \left[A_{n1}\right] & \ldots & \left[A_{nn}\right]\end{array}\right) \; ,
\end{eqnarray}
where every block is a tridiagonal matrix.
The diagonal blocks have the following structure
\begin{eqnarray} \label{e17} \hspace{-2ex}
 [A_{ii}] = \left( \!\begin{array}{ccccccc}
  1\!+\!\frac{2}{h^2}a_{ii} & -\frac{1}{h^2}a_{ii} & 0 & 0& 0& \ldots & 0 \\
  -\frac{1}{h^2}a_{ii} & 1\!+\!\frac{2}{h^2}a_{ii} & -\frac{1}{h^2}a_{ii} & 0 & 0 & \ldots & 0 \\
  0 & -\frac{1}{h^2}a_{ii} & 1\!+\!\frac{2}{h^2}a_{ii} & -\frac{1}{h^2}a_{ii} & 0 & \ldots & 0 \\
  \vdots &&&&&& \vdots \\
  0 & 0 & 0 & \ldots & -\frac{1}{h^2}a_{ii} & 1\!+\!\frac{2}{h^2}a_{ii} & -\frac{1}{h^2}a_{ii}  \\
  0 & 0 & 0 & 0 & \ldots & -\frac{1}{h^2}a_{ii} & 1\!+\!\frac{2}{h^2}a_{ii} 
\end{array}\!\right),
\end{eqnarray}
and nondiagonal ones
\begin{eqnarray} \label{e18}
 [A_{ij}] = \left( \begin{array}{ccccccc}
  \frac{2}{h^2}a_{ij} & -\frac{1}{h^2}a_{ij} & 0 & 0& 0& \ldots & 0 \\
  -\frac{1}{h^2}a_{ij} & \frac{2}{h^2}a_{ij} & -\frac{1}{h^2}a_{ij} & 0 & 0 & \ldots & 0 \\
  0 & -\frac{1}{h^2}a_{ij} & \frac{2}{h^2}a_{ij} & -\frac{1}{h^2}a_{ij} & 0 & \ldots & 0 \\
  \vdots &&&&&& \vdots \\
  0 & 0 & 0 & \ldots & -\frac{1}{h^2}a_{ij} & \frac{2}{h^2}a_{ij} & -\frac{1}{h^2}a_{ij}  \\
  0 & 0 & 0 & 0 & \ldots & -\frac{1}{h^2}a_{ij} & \frac{2}{h^2}a_{ij} 
\end{array}\right).
\end{eqnarray}
The set of linear equations (\ref{e14}) can be solved by standard methods,
for instance the LU decomposition \cite{PTVF}.

\section{Stochastic integral}

The essential part of the mild solution (\ref{eq3})
is the stochastic convolution (\ref{eq4}). 
We present here a particular case when the resolvent $S(t)$ is known
in analytical form. 
Let us focus the attention on stochastic Volterra equation (\ref{eq2})
with the function $a$ in the form $a(t)=t^{\alpha-1}/\Gamma(\alpha)$.
This is the integrodifferential equation \cite{Fu},\cite{Ka1},\cite{Pr}.
For three particular cases, $\alpha=0,1,2$ the analytical form
of the resolvent $S$ is known:
\begin{equation}\label{eq19}
 (S(t)v)(x) := \int_{-\infty}^{\infty} \phi_\alpha (t,x-y)v(y)dy
  = \int_{-\infty}^{\infty} \phi_\alpha (t,y)v(x-y)dy,
\end{equation}
where the last form comes from the convolution property.
We shall illustrate the applicability of our numerical algorithms
with two cases of the above $a$ functions, the case with  $\alpha=1$
and with  $\alpha=2$. Then the function $\phi_\alpha$ in (\ref{eq19})
takes the following form
\begin{eqnarray}\label{eq20}
 \phi_1 (t,x) & = & \frac{1}{\sqrt{4\pi t}}\,\exp(-\frac{x^2}{4t}) 
 \hspace{15.2ex} \mbox{for}~~\alpha=1 \\
 \phi_2 (t,x) & = & \half(\delta(t-x)+\delta(t+x)) 
 \hspace{10ex} \mbox{for}~~\alpha=2 . 
\end{eqnarray}

We assume the process in the form $W(t,x)=W_1(t)\,W_2(x)$.
Then, the algorithm for an approximate construction of the stochastic 
integral (\ref{eq4}) can be built in the following way.
Let us introduce a time grid $\{t_i=i\tau:~i=0,1,\ldots,I\}$ on
$[0,T]$, i.e.\ $\tau=T/I$ and next a finite sequence of independent 
random variables $\{\zeta_i\},~i=1,2,\ldots,I$, with standard normal
distribution.
The approximation for the convolution (\ref{eq4}) can be written in the
form
 \begin{eqnarray}\label{eq21}
 \int_0^t S(t-s)\,dW(s,x) & = & \sum_{i=0}^{I-1} S(t-s_i)
 [W(s_{i+1},x)-W(s_{i},x)]  \\
 & = & \tau^{\half} \sum_{i=0}^{I-1} \zeta_i 
 \int_{-\infty}^{\infty} \phi_\alpha (t-s_i,x-y)\,W_2(y) dy
 \nonumber
\end{eqnarray}
For further specification we choose $W_2(x)=C\,X_0(x)=C\,e^{-x^2/4}$
(the constant $C$ represents a 'strength' of stochastic forces).
With this assumption, after performing the integral (\ref{eq21})
for particular $\phi_\alpha$ one obtains 
 \begin{eqnarray}\label{eq22}
 \lefteqn{\int_0^t S(t-s)\,dW(s,x) =} \\
 && \hspace{10ex} C\tau^{\half} \sum_{i=0}^{I-1} \zeta_i 
 \frac{1}{\sqrt{1+t+s_i}}\,\exp (\frac{-x^2}{4(1+t+s_i)}) ~~~
 \mbox{for}~~~\alpha=1,\nonumber
\end{eqnarray}
and
\vfill
\newpage
 \begin{eqnarray}\label{eq23}
 \lefteqn{\int_0^t S(t-s)\,dW(s,x) =}  \\
 & &~ C\tau^{\half} \sum_{i=0}^{I-1} \zeta_i  \half \left[
 \exp (\frac{-(x-t+s_i)^2}{4})
 +\exp (\frac{-(x+t-s_i)^2}{4})\right] ~~~
 \mbox{for}~~~\alpha=2.\nonumber
\end{eqnarray}
These explicite forms were inserted into the numerical code.
The sequence of independent random variables 
$\{\zeta_i\},~i=1,2,\ldots,I$, with standard normal distribution
was generated using subroutines 
{\tt gasdev} and {\tt ran1} from \cite{PTVF}. 

\section{Numerical results}

We illustrate the efficiency of the numerical approach on two 
examples of the function $a$, mentioned earlier.
As the initial value of the $X(t=0,x)=X_0(x)$ we chose the Gaussian
$X_0(x)= e^{-x^2/4}$. The grid in $x$ variable contained $m=150$ 
intervals with $h=0.2$, covering the interval $x\in[-15,15]$.  
The dimension of the approximation subspace in the Galerkin method
was chosen as $n=8$. The resulting  dimension of the matrix
$\mathcal{A}$ was then 1208$\times$1208 and
calculations were performed up to $T=6$.

\begin{figure}[bht]  
\begin{center}
\resizebox{0.95\textwidth}{!}{\includegraphics{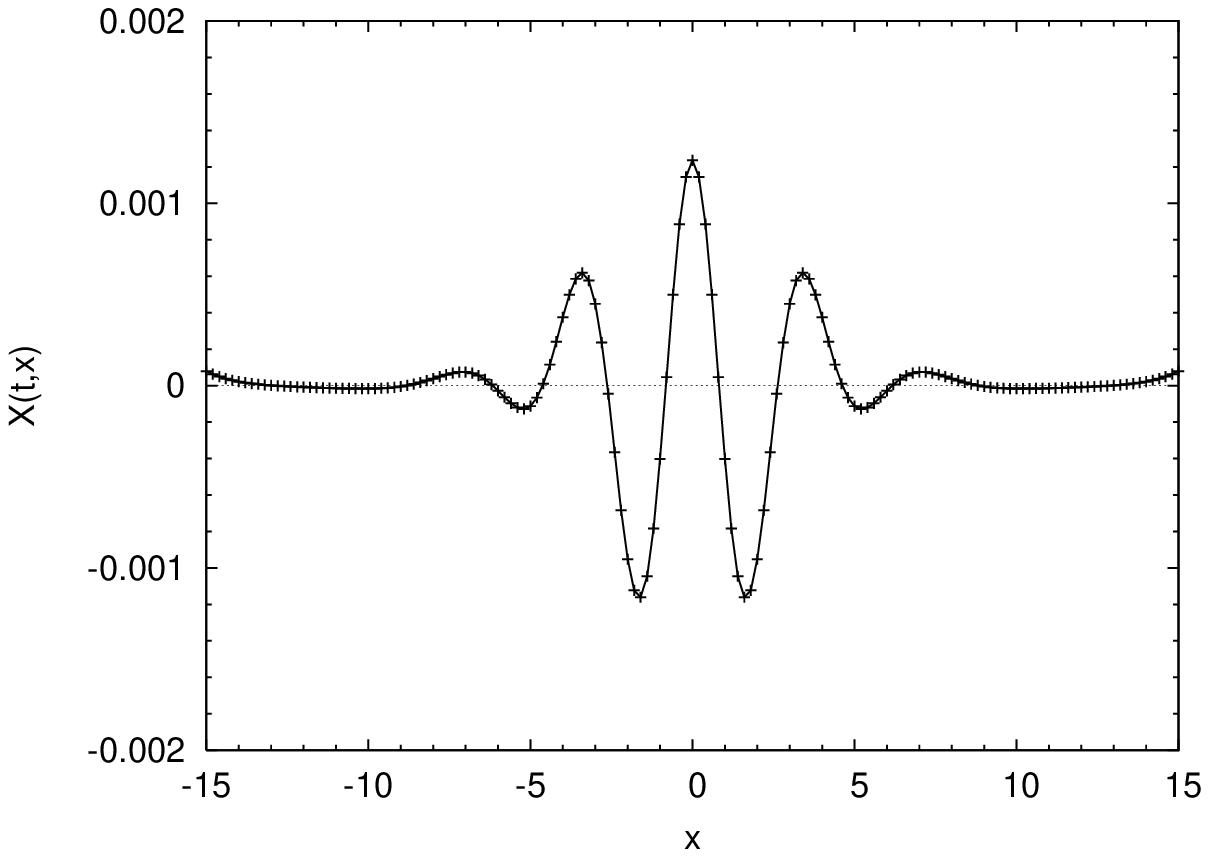}}
\resizebox{0.95\textwidth}{!}{\includegraphics{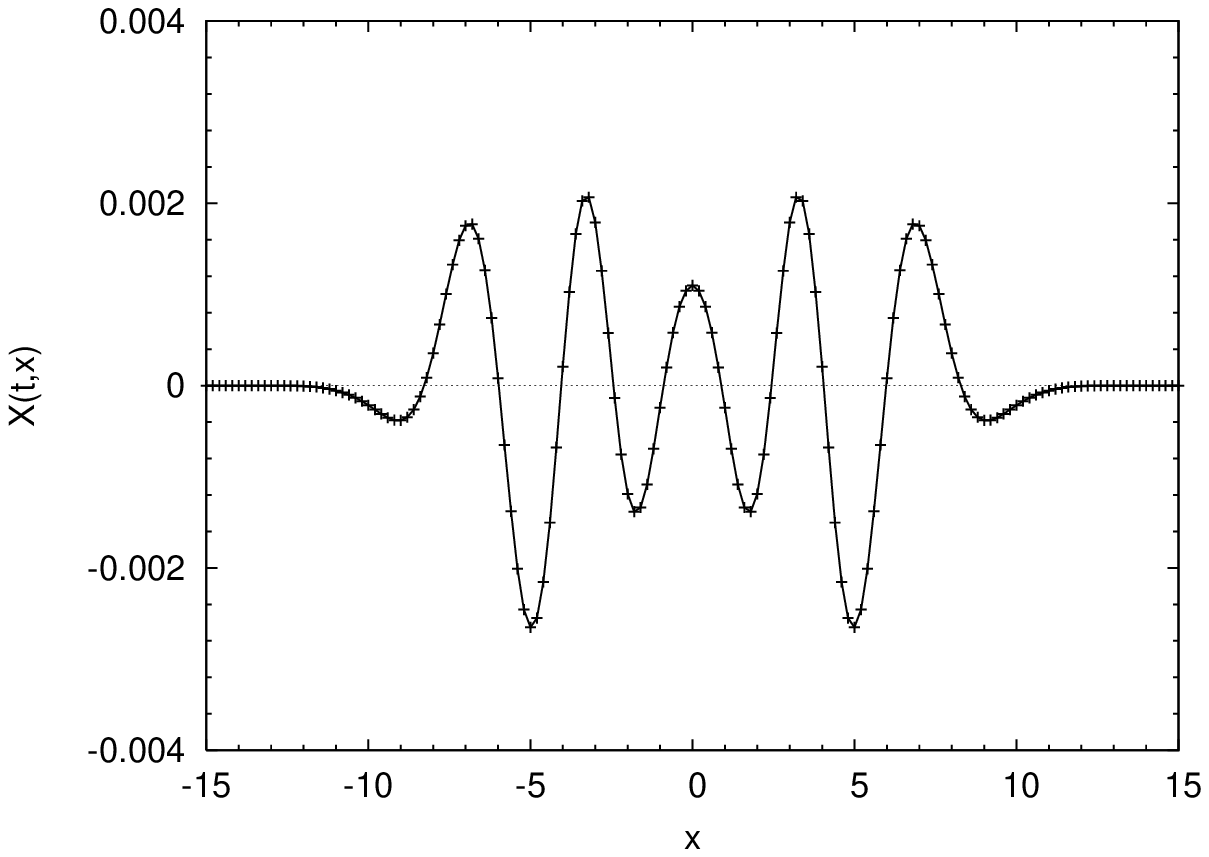}}
\end{center}
\caption{Errors of the numerical solutions to deterministic Volterra
equation (\ref{eq1a}) at $t=6$. Top: $\alpha=1$ case, bottom: $\alpha=2$ case.
The symbols represent the difference between the exact (analytical)
solution and the numerical solution at given grid points.}
\label{err}
\end{figure}

In fig.~\ref{err} we show the errors of the numerical solutions  
to the deterministic equation (\ref{eq1a}) obtained in cases
$\alpha=1$ (top), and  $\alpha=2$ (bottom). In both cases the
errors are relatively small.

\begin{figure}[hbt]  
\begin{center}
\resizebox{0.95\textwidth}{!}{\includegraphics{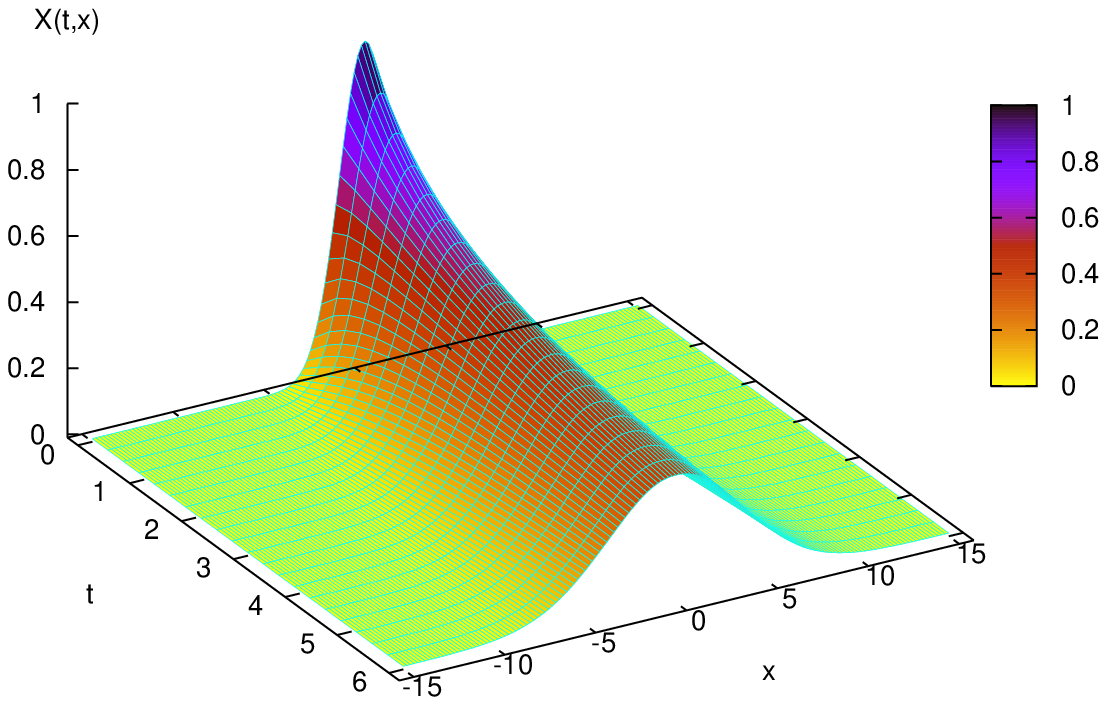}}
\resizebox{0.95\textwidth}{!}{\includegraphics{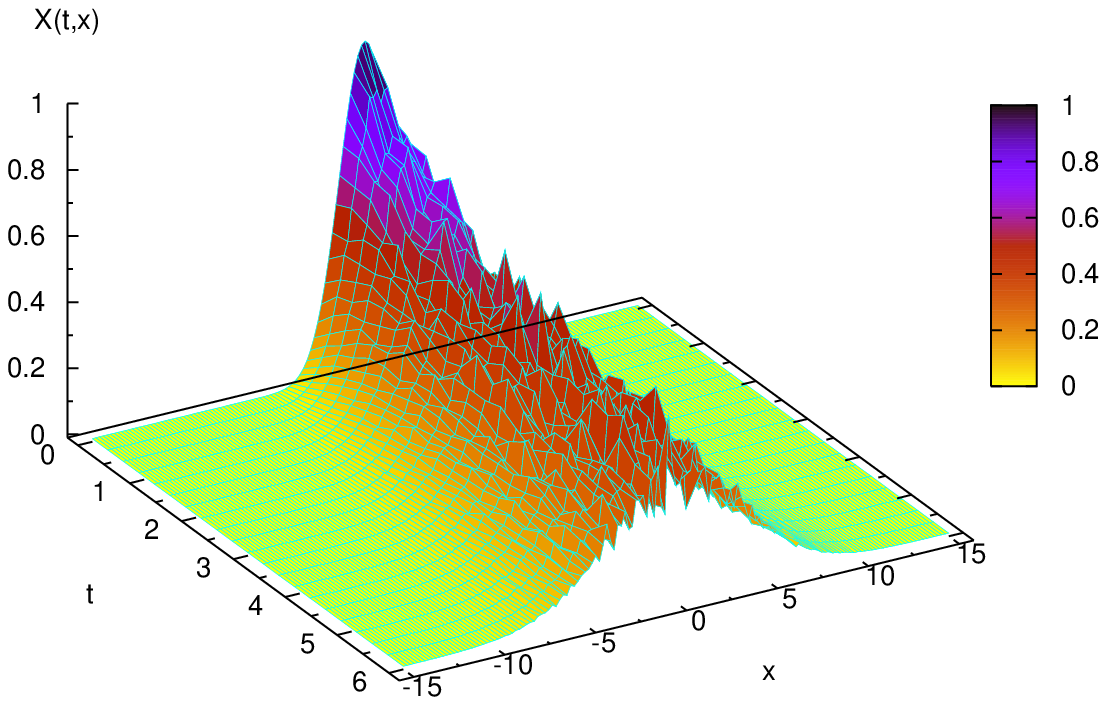}}
\end{center}
\caption{Numerical solution to  Volterra equation with 
$\alpha=1$ for $t\in [0,6]$:  
the deterministic solution (top) and the stochastic one (bottom).}
\label{a1ds}
\end{figure}

The top part of the fig.~\ref{a1ds} displays the solution of the 
deterministic Volterra equation \ref{eq1} for the case $\alpha=1$
as function of time, $t\in[0,6]$. The bottom part of the  
fig.~\ref{a1ds} shows the example of a single stochastic trajectory
(i.e.\ the sum of deterministic solution and stochastic integral)
for the same case. 

\begin{figure}[hbt] 
\begin{center}
\resizebox{0.95\textwidth}{!}{\includegraphics{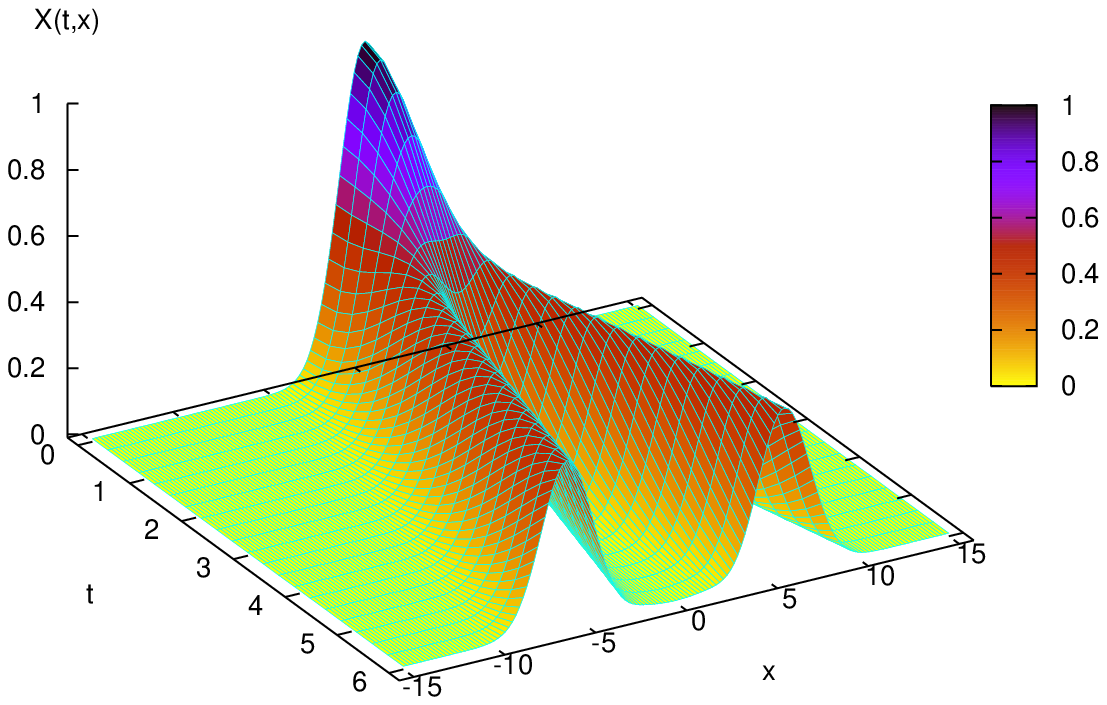}}
\resizebox{0.95\textwidth}{!}{\includegraphics{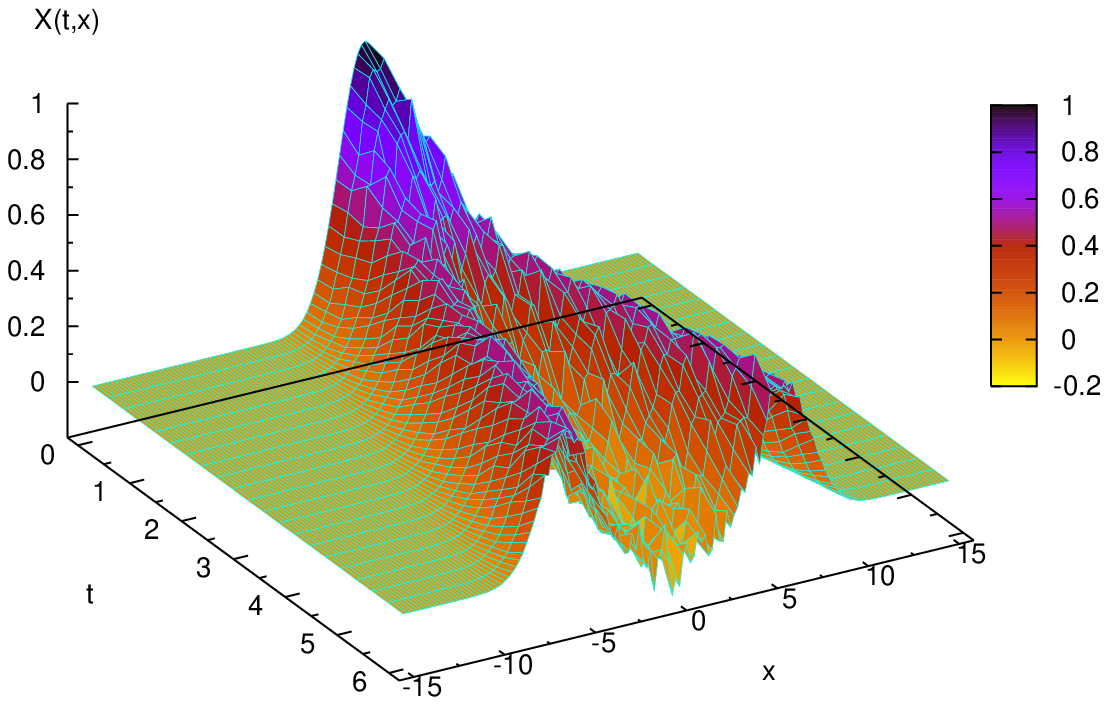}}
\end{center}
\caption{Numerical solution to  Volterra equation with 
$\alpha=2$ for $t\in [0,6]$:  
the deterministic solution (top) and the stochastic one (bottom).}
\label{a2ds}
\end{figure}

The fig.~\ref{a2ds} present the corresponding
solutions for the case $\alpha=2$.
For stochastic convolutions the value of the constant 
$C$ was chosen to be $C=0.1$.

\end{document}